\documentstyle[12pt,leqno,epsfig]{article}
\renewcommand{\arraystretch}{1.0}

\newfont{\myfontsmall}{msbm10 scaled 1200}

\def\bpn{\bigskip\par\noindent}
\def\mpn{\medskip\par\noindent}
\def\bea{\begin{eqnarray}}
\def\eea{\end{eqnarray}}
\def\beas{\begin{eqnarray*}}
\def\eeas{\end{eqnarray*}}
\def\be{\begin{equation}}
\def\ee{\end{equation}}
\def\TN{T_{N}}
\def\TA{T_{N}^{*}}

\def\supp{\mbox{supp\/}}
\def\.{{\;}}
\def\BN{{\cal B}_N}
\def\BM{{\cal B}_M}

\def\z{\mbox{\myfontsmall Z}}
\def\r{\mbox{\myfontsmall R}}

\def\n{\mbox{\myfontsmall N}}

\begin{document}
\newtheorem{theorem}{Theorem}[section]
\newtheorem{proposition}[theorem]{Proposition}
\newtheorem{algorithm}[theorem]{Algorithm}
\newtheorem{definition}[theorem]{Definition}
\newtheorem{corollary}[theorem]{Corollary}
\newtheorem{remark}[theorem]{Remark}
\newtheorem{lemma}[theorem]{Lemma}
\newtheorem{example}[theorem]{Example}
\newtheorem{assumption}[theorem]{Assumption}

\renewcommand{\theequation}{\thesection.\arabic{equation}}
\def\abschnitt{\setcounter{equation}{0} \section}
\renewcommand{\thefigure}{\thesection.\arabic{figure}}
\setcounter{totalnumber}{4}

\newcounter{oldasseqn}

\newcounter{saveeqn}
\newcommand{\asseqn}{\setcounter{saveeqn}{\value{equation}}
 \setcounter{equation}{\value{oldasseqn}}
 \renewcommand{\theequation}{\mbox{P\arabic{equation}}}}

\newcommand{\reseteqn}{\setcounter{oldasseqn}{\value{equation}}
 \setcounter{equation}{\value{saveeqn}}
 \renewcommand{\theequation}{\thesection.\arabic{equation}}}

\vspace*{-2cm}
\noindent
{\small \tt Published in Num.Funct.Anal.Opt., 19/3-4: 353--375, 1998.}

\vspace*{1cm}

\Large
{\bf
\begin{center}
\vspace*{10mm}
A MULTI--LEVEL ALGORITHM FOR THE SOLUTION OF MOMENT PROBLEMS
\end{center}
}
\normalsize
\vspace{0.3cm}
\begin{center}
{
OTMAR SCHERZER\\
Institut f\"ur Industriemathematik, Universit\"at Linz,\\
Altenberger Str. 69, A--4040 Linz, Austria\\
\bigskip
THOMAS STROHMER\\
Institut f\"ur Mathematik, Universit\"at Wien,\\
Strudlhofgasse 4, A-1090 Wien, Austria }
\end{center}

\begin{abstract}
We study numerical methods for the solution of general linear 
moment problems, where the solution belongs to a family of nested 
subspaces of a Hilbert space. Multi-level algorithms, based on 
the conjugate gradient method and the Landweber--Richardson method 
are proposed that determine the ``optimal'' reconstruction level
a posteriori from quantities that arise during the numerical 
calculations. As an important example we discuss the reconstruction
of band-limited signals from irregularly spaced noisy samples, 
when the actual bandwidth of the signal is not available. Numerical 
examples show the usefulness of the proposed algorithms. 
\end{abstract}

\mpn


\abschnitt{Introduction}
\label{sec1}

We study numerical methods for the solution of a family of
general linear moment problems
\begin{equation}
\label{1.1}
\langle x,g^N_j \rangle = \mu_j \qquad j \in I
\end{equation}
where 
$\{g_j^N\}_{j \in I}$ is a complete system in $X_N$ and
$\{X_N\}_{j \in I}$ is a family of nested subspaces of a separable
Hilbert space $X$. $\langle x,g_j^N \rangle$ denotes the 
$j$--th moment. We consider the problem of recovering $x$ 
from noisy observations of the moments $\mu_j +\delta_j$.
\bpn
Most existing methods for moments problems require exact
a priori information about the space $X_N$ to which $x$ 
belongs to \cite{BacGil68,XiaNas94,LM91}. In many applications 
however this knowledge is not available, often we only know that
$x$ belongs to some space of a whole family of 
Hilbert spaces $X_N$. 

An important example is the
problem of reconstructing a band-limited signal $x$ from 
irregularly spaced, noisy samples $x(t_j) +\delta_j$, without
explicit knowledge on the bandwidth of $x$. In other words we know 
that $x$ belongs to a space of a family of band-limited 
functions $X_N$ (where $N$ represents the bandwidth), but we do 
not know to which one. Most reconstruction algorithms require the 
a priori information of the actual bandwidth of $x$ (or at 
least a very good estimate) to be applicable. In numerical
reconstructions overestimating the bandwidth of $x$ may result in a
highly oscillating function, due to the noise in the samples.
Underestimating the bandwidth may give a poor approximation
to the original signal. From this point of view a 
proper choice of the level $N$ can also be seen as a 
regularization procedure.
\bpn
For the numerical solution of the family of general moment
problems we propose a multi--level approach that determines
an ``optimal'' level $N$ a--posteriori from
quantities that arise during the numerical calculations.
The proposed algorithm does not use a--priori information on
the level $N$ of the solution -- all information
required
in the analysis of this algorithm is that there exists a level
(bandwidth) $N_*$ such that the solution $x$ of (\ref{1.1})
satisfies $x_* \in X_{N_*}$. 
\bpn
The outline of this paper is the following: In Section
\ref{sec2}
and Section \ref{sec3} we study multi--level algorithms for the
efficient stable solution of moment problems. In particular
we concentrate on a multi-level conjugate gradient algorithm
and a multi-level Landweber--Richardson algorithm.
These algorithms are formulated in a general setting.
In Section \ref{sec4} several applications
to the numerical solution of moment problems are given,
with the emphasize on reproducing kernel Hilbert spaces. In
Section~\ref{sec5} we demonstrate in detail how the proposed 
algorithms can be applied to the problem of reconstructing a 
band-limited signal of
unknown bandwidth from noisy irregularly spaced samples.
Finally in Section~\ref{sec6} we demonstrate the performance 
of the multi-level algorithms by some numerical examples.

\abschnitt{A multi--level algorithm based on the conjugate gradient method}
\label{sec2}
In this section we introduce a multi--level algorithm based on the 
conjugate gradient for the numerical solution of moment problems. 
In this algorithm the levels are adapted inductively and at each 
level a conjugate gradient method is implemented. 
\bpn
In order to make the paper self--contained we outline the conjugate 
gradient method for the solution of linear systems: Let $T : X \to Y$ 
be a linear operator between separable Hilbert spaces $X$ and $Y$. 
By $T^*$ we denote the adjoint of the operator $T$. One efficient method 
for the solution of $Tx=y$ is the conjugate gradient method applied to 
the normal equation $T^* x T = T^* y$. This method, which is denoted by 
CGNE in the literature, is defined as follows 
(see e.g. \cite{Han95,Han96,EngHanNeu96})
\beas
\label{CGNE}
\fbox{$
\begin{array}{l}
\begin{array}{rcl}
x_0 &=& 0\\
r_0 &=& y\\
w_0 &=& y\\
k        &=& 0
\end{array}
\\
\mbox{ repeat }\\
\left\{
\begin{array}{rcl}
d_k &=& T^* w_k\\
\alpha_k &=& \|T^* r_k\|^2/\|T d_k\|^2\\
x_{k+1} &=& x_k + \alpha_k d_k\\
r_{k+1} &=& r_k - \alpha_k T d_k\\
\beta_k &=&  \|T^* r_{k+1}\|^2/\|T^* r_k\|^2\\
w_{k+1} &=& r_{k+1} + \beta_k w_k\\
k&=&k+1\\
\end{array}
\right.
\\
\mbox{ until } \|r_k\| < \rho \|y\|
\end{array}
$}
\eeas
\bpn
Before we study a multi--level scheme based on the CGNE method
we introduce some notation:
\begin{assumption}
\label{ass2.1CG}
$\{X_N\}_{N \in \n_0}$ is a family of linearly ordered subspaces of a 
separable Hilbert space $X$. For $N \in \n_0$ let $P_N$ be the 
orthogonal projectors from $X$ into $X_N$ and let ${\cal R}(\TA)$
denote the range of $\TA$. For $N \in \n_0$ the family of linear operators 
$\TN: X_N \to Y$ satisfies
\bea
\label{2.1CG}
\left\{
\begin{array}{cl}
\mbox{\em (i)}   & \,\,\,\|\TN\| \leq 1\,,\\
\mbox{\em (ii)}  & \,\,\,\TN x_M = T_M x_M \mbox{ for all } N \geq M\,,\\
\mbox{\em (iii)} & \,\,\,{\cal{R}}(\TA) \subseteq X_N\,.
\end{array}
\right.
\eea
\end{assumption}
The first assumption in (\ref{2.1CG}) is that the operators $\TN$ 
are properly scaled. 
The second assumption is that on subspaces of level $M (< N)$ the 
operators $T_M$ and $\TN$ coincide. The third assumptions guarantees 
that the iterates of the CGNE method are contained in $X_N$. 
\bpn
For given data $y \in Y$ we are concerned with finding $N \in \n_0$ and 
$x \in X_N$, which satisfy
\be
\label{2.2CG}
\TN x = y\..
\ee
The terminology of a solution of (\ref{2.2CG}) is well--defined by 
Assumption (\ref{ass2.1CG}), since from this it follows that a solution 
on a lower scale is also a solution on a higher scale. 
\bpn
In the sequel we will assume that there exists a solution of (\ref{2.2CG}).
\begin{assumption}
\label{ass2.2CG}
There exists $N_* \in \n_0$ and $x_* \in X_{N_*}$ which satisfy 
(\ref{2.2CG}). 
\end{assumption}
We note that we do not assume the knowledge of $N_*$. 
If for instance $X_{N}$ denotes the Payley--Wiener space of band-limited 
functions of bandwidth less or equal $N$ and $N \to \infty$, then 
the meaning of Assumption \ref{ass2.2CG} is: There exists a band-limited 
function which solves (\ref{2.2CG}). 
\bpn
By $y^\delta \in Y$ we denote 
measured data of $y$, which satisfies 
\be
\label{2.3CG}
\|y^\delta - y\| \leq \delta\..
\ee
\begin{algorithm}
\label{alg2.3CG} 
Let $\eta > 0$, $\tau > 1$.
\begin{description}
\item{At level $0:$} Without loss of generality we assume 
      that $x_{0,0}^\delta \in X_{0}$ satisfies 
      \begin{enumerate}
      \item \label{2-1CG}
            $\|T_0 (x_{0,0}^\delta) - y^\delta\|^2 > 
            2(1+\eta)\left( \delta + \|P_{0}x_* - x_*\| \right)\|y^\delta\|$,
      \item \label{2-2CG}
            $\|T_0 (x_{0,0}^\delta) - y^\delta\|^2 > 
            2(1+\eta) \tau \delta \|y^\delta\|$.
      \end{enumerate}
      If \ref{2-2CG}. does not hold, then $x_{0,0}^\delta$ is 
      accepted as an approximate solution and the multi--level 
      algorithm is terminated. 
      If \ref{2-1CG}. does not hold, but \ref{2-2CG}. is satisfied, 
      then we can choose a level $J$ such that 
      $$\|T_J (x_{0,0}^\delta) - y^\delta\|^2 =
        \|T_0 (x_{0,0}^\delta) - y^\delta\|^2 > 
         2(1+\eta)\left( \delta + \|P_Jx_* - x_*\| \right)\|y^\delta\|\..$$
      In this situation we simply renumber the levels and make the 
      setting $0=J$. 
      For the new setting $J=0$ both \ref{2-1CG}. and \ref{2-2CG}. hold, 
      and the CGNE method at level $0$ is performed.
      The iteration at scale $X_0$ is terminated if for the first 
      time
      $$
      \|T_0 (x_{k+1,0}^\delta) - y^\delta\|^2 \leq  
        2(1+\eta)\left( \delta + \|P_0x_* - x_*\| \right)\|y^\delta\|\..
      $$
      The termination index is denoted by $k=k_*(0,\delta)$.
\item{From level $N$ to level $N+1$:} 
      If $\|T_{N+1} (x_{k_*(N,\delta)+1,N}^\delta) - y^\delta\|^2 \leq 
          2(1+\eta) \tau \delta \|y^\delta\|$, 
      then the multi--level algorithm terminates and 
      $x_{k_*(N,\delta)+1,N}^\delta$ is an 
      approximate solution.
      Else we define $x_{k_*(N,\delta),N}^\delta =: x_{0,N+1}^\delta$.
      Since $\|P_{N+1}x_* - x_*\| \leq \|P_Nx_* - x_*\|$ it follows 
      from the definition of $k_*(N,\delta)$:
      $$\|T_{N+1} (x_{0,N+1}^\delta) - y^\delta\| > 
        2(1+\eta)\left( \delta + \|P_{N+1}x_* - x_*\| \right)\|y^\delta\|
        \..$$
      Therefore, at least one iteration of CGNE at level $N+1$ is 
      performed.
      The iteration at scale $X_{N+1}$ is terminated if 
      $$\|T_{N+1} (x_{k+1,N+1}^\delta) - y^\delta\|^2 \leq 
        2(1+\eta)\left( \delta + \|P_{N+1}x_* - x_*\| \right)\|y^\delta\|$$
      for the first time. The termination index 
      is denoted by $k=k_*(N+1,\delta)$.
\end{description}
\end{algorithm}

The idea of using multi--level algorithms as considered in this 
paper for the solution of ill--posed problems has been considered 
in \cite{Sch96}.

In the following, we verify some basic monotonicity properties of the 
iterates of the multi--level algorithm \ref{alg2.3CG}. Note
that in practice when the CGNE--method is used for solving linear systems
the CGNE method is usually terminated with a discrepancy principle 
(see e.g. \cite{Han96}). Hanke \cite{Han96} showed by a counterexample 
that in general the error of the iterates of the CGNE method
does not monotonically decrease when the 
iteration is terminated by the discrepancy principle. 
In this paper the CGNE--method is combined with a multi--level algorithm 
and we expect that a reasonable performance of this method is depending 
on providing a good initially guess for starting the CGNE method at a 
higher level. Therefore, we concentrate in this paper on a stopping 
criterion, which guarantees monotonicity of the iterates at each level.
\begin{proposition}
\label{pr2.4CG}
Let assumption \ref{ass2.1CG} hold, and let the iteration at fixed level $N$ 
be stopped according to the generalized discrepancy principle 
\be
\label{2.4CG}
\| \TN (x_{k,N}^\delta) - y^\delta\|^2 > 
2(1+\eta) \left( \delta + \|P_Nx_*-x_*\| \right)\|y^\delta\|\, \mbox{ for }
0 \leq k \leq k_*(N,\delta)\..
\ee
Then $k_*(N,\delta) < \infty$ and 
\beas
\begin{array}{rclr}
\|x_{j+1}^N - P_Nx_*\|  &\leq &
\|x_{j}^N - P_Nx_*\| & \quad j=0,...,k_*(N,\delta)\,,\\
\|x_{j+1}^N - x_*\| & \leq &
\|x_{j}^N - x_*\| & \quad j=0,...,k_*(N,\delta)\;.
\end{array}
\eeas
\end{proposition}
{\bf{Proof:}} Application of the general results in \cite{Han96} 
immediately yields $k_*(N,\delta) < \infty$ and for 
$j=0,...,k_*(N,\delta)$
$\|x_{j+1}^N - x_*\| \leq \|x_{j}^N - x_*\|$. 
Therefore we have 
$$\|x_{j+1}^N - P^N x_*\|^2 + 2\langle x_{j+1}^N - P^N x_*,P^N x_*-x_* \rangle
  \leq \|x_{j}^N - x_*\|^2 + \langle x_{j+1}^N - P^N x_*,P^N x_*-x_* \rangle.$$ 
From the remark following Assumption \ref{ass2.1CG} it follows 
$x_{j}^N, x_{j+1}^N \in X_N$ and therefore the assertion is proven.
\hfill $\Box$ 
\bpn
From Proposition \ref{pr2.4CG} it follows that (\ref{2.4CG}) determines a 
well--defined stopping index for the CGNE method at level $N$.
\subsection{\label{sec2.1CG} Convergence and stability of the multi--level 
                           algorithm}
Using the monotonicity results of the iterates of the multi--level schemes 
we prove below that the multi--level scheme is convergent and stable.
The following lemma will be central for our further considerations.
\begin{lemma}
\label{le2.5CG}
Let $\delta = 0$. If Assumptions \ref{ass2.1CG} and \ref{ass2.2CG} 
holds, then the multi--level algorithm \ref{alg2.3CG}
\begin{description}
\item{ } terminates after a finite number of iterations
\item{ or }
\item{ } 
      there exists a level $\hat{N}$ such that at each level
      $N < \hat{N}$ a finite number of CGNE--iterations
      are performed and at level $\hat{N}$ the 
      CGNE iteration does not terminate.
\end{description} 
\end{lemma}
{\bf{Proof:}}
If $\delta = 0$, then the stopping criterion (\ref{2.4CG}) 
for the CGNE method at level $N$ is
\be
\label{2.15CG}
\|\TN(x_{k,N}^\delta) - y\|^2 > 
2(1+\eta)\| P_N x_* - x_*\|\|y\| \mbox{ for } 
0 \leq k \leq  k_*(N,0)\..
\ee
If $\delta = 0$, then the stopping criterion for the multi--level 
algorithm is
\be
\label{2.16CG}
\|\TN(x_{k,N}^\delta) - y\| = 0\..
\ee
Let $\hat{N}$ be the first level--index where $P_Nx_* = x_*$. A finite 
level--index exists by Assumption \ref{ass2.2CG}.
From Proposition \ref{pr2.4CG} it follows that for each level $N < \hat{N}$ 
the number of iterations of the CGNE method is finite, 
which yields the assertion. \hfill $\Box$
\bpn
Using the previous lemma, we are able to prove that the multi--level 
CGNE algorithm is convergent in the case of exact 
data.
\begin{theorem}
\label{th2.6CG}
Let Assumption \ref{ass2.2CG} hold. Then the multi--level CGNE 
algorithm converges to a solution $\hat{x}$. 
\end{theorem}
{\bf{Proof:}}
Let $\hat{N}$ be the smallest index $N$ such that $P_Nx_* = x_*$.
We differ between two cases:
\begin{enumerate}
\item The multi--level CGNE method terminates at level 
      $\hat{N}$ after a finite number of iterations
\item The multi--level CGNE method does not terminate
\end{enumerate}
In the first case it follows from (\ref{2.16CG}) that the last 
iterate $\hat{x}$ satisfies $ T_{\hat{N}}^* \hat{x} = y\;,$
which shows that $\hat{x}$ is the solution.\\
In the second case it follows immediately from Theorem 3.4 in 
\cite{Han95} that $x_n$ is convergent to a solution.
\hfill$\Box$
\bpn
In the case of perturbed data, it follows from Proposition \ref{pr2.4CG} 
that the discrepancy principle (\ref{2.4CG}) determines a well--defined 
(finite) stopping index at level $N < \hat{N}$, which is denoted by 
$k_*(N,\delta)$. 
In the following, we will assume that the multi--level algorithm is 
finally terminated if for the first time
\be
\label{2.7CG}
\|\TN(x_{k,N}^\delta) - y^\delta\|^2 \leq 
        2(1+\eta)\tau\delta\|y^\delta\|\..
\ee
From Lemma \ref{le2.5CG} and Assumption \ref{ass2.2CG} it follows that 
this criterion determines a well--defined stopping index at level 
$N(\delta)$, which is denoted by $K_*(N(\delta),\delta)$.
\bpn
Our next result shows that the multi--level iteration is a stable method.
\begin{theorem}
\label{th2.7CG}
Let the assumptions of Theorem \ref{th2.6CG} hold and assume that 
$T_{N_*}$ is continuously invertible (here $N_*$ is as
introduced in Assumption \ref{ass2.2CG}).
If $y^\delta$ fulfills (\ref{2.3CG}), and if the iteration at level $N$ 
is stopped according to the discrepancy principle (\ref{2.4CG}), and the 
multi--level iteration is terminated by (\ref{2.7CG}), then 
$$x_{N_n,K_*(N_n,\delta)}^\delta \to x_*, \qquad \delta \to 0\..$$
\end{theorem}  
{\bf{Proof:}} Let $\delta_n$, $n=1,2,...,$ be a sequence converging to 
zero as $n \to \infty$, and let $y_n:=y^{\delta_n}$ be a corresponding 
sequence of
perturbed data. For each pair $(\delta_n,y_n)$, denote by 
$K_n=K_*(N(\delta_n),\delta_n)$ the corresponding stopping index 
at termination level $N_n = N(\delta_n)$ determined by (\ref{2.7CG}).
A comparison of $\delta_n + \|P_{N(\delta_n)}x_* - x_*\|$ and 
$\tau\delta_n$ shows that $N(\delta_n) \leq N_*$. 
From Proposition \ref{pr2.4CG} and the results in \cite{Han96} it follows 
that for 
$0 \leq k \leq k_*(N,\delta)$ and $0 \leq N \leq N_*$
$$\|x_{k,N}^\delta - x_*\| \leq \|x_{0,0}^\delta - x_*\|
  \mbox{ and } 
  \|T_{N_*}x_{k,N}^\delta - y\| \to 0\,.$$
From the continuous invertibility of $T_{N_*}$ the assertion 
follows. \hfill$\Box$
\bpn

\abschnitt{A multi--level algorithm based on the Landweber--Richardson 
           method}
\label{sec3}
As an alternative to the multi--level algorithm based on CGNE we study in this
section a multi--level iteration based on the Landweber--Richardson 
method for solving moment problems. It will turn out that different 
stopping criteria can be used for stopping the Landweber--Richardson 
method than for the CGNE method.
\bpn
The multi--level iteration based on the Landweber--Richardson method 
is defined as follows
\begin{algorithm}
\label{alg3.1} 
Let $\eta > 0$, $\tau > 1$.
\begin{description}
\item{At level $0:$} Without loss of generality we assume 
      that $x_{0,0}^\delta \in X_{0}$ satisfies 
      \begin{enumerate}
      \item \label{3-1}
            $\|T_0 (x_{0,0}^\delta) - y^\delta\| > 
            2(1+\eta)\left( \delta + \|P_{0}x_* - x_*\| \right)$,
      \item \label{3-2}
            $\|T_0 (x_{0,0}^\delta) - y^\delta\| > 
            2(1+\eta) \tau \delta$.
      \end{enumerate}
      If \ref{3-2}. does not hold, then $x_{0,0}^\delta$ is 
      accepted as an approximate solution and the multi--level 
      algorithm is terminated. 
      If \ref{3-1}. does not hold, but \ref{3-2}. is satisfied, 
      then we can choose a level $J$ such that 
      $$\|T_J (x_{0,0}^\delta) - y^\delta\| =
        \|T_0 (x_{0,0}^\delta) - y^\delta\| > 
         2(1+\eta)\left( \delta + \|P_Jx_* - x_*\| \right)\..$$
      In this situation we simply renumber the levels and set $0=J$. 
      With the new setting both \ref{3-1}. and \ref{3-2}. hold, 
      and the ``Landweber--Richardson'' iteration at level $0$ is 
      performed
      $$x_{k+1,0}^\delta = x_{k,0}^\delta - 
        T_0^* \left( T_0 (x_{k,0}^\delta)-y^\delta \right)\..$$
      The iteration at scale $X_0$ is terminated if for the first 
      time
      $$
      \|T_0(x_{k+1,0}^\delta) - y^\delta\| \leq  
        2(1+\eta)\left( \delta + \|P_0x_* - x_*\| \right)\..
      $$
      The termination index is denoted by $l_*(0,\delta)$.
\item{From level $N$ to level $N+1$:} 
      If $\|T_{N+1} (x_{l_*(N,\delta)+1,N}^\delta) - y^\delta\| \leq 
          2(1+\eta) \tau \delta$, 
      then the multi--level algorithm terminates and 
      $x_{l_*(N,\delta)+1,N}^\delta$ is an 
      approximate solution.
      Else we define $x_{l_*(N,\delta),N}^\delta =: x_{0,N+1}^\delta$.
      Since $\|P_{N+1}x_* - x_*\| \leq \|P_Nx_* - x_*\|$ it follows from the 
      definition of $l_*(N,\delta)$:
      $$\|T_{N+1}(x_{0,N+1}^\delta) - y^\delta\| > 
        2(1+\eta)\left( \delta + \|P_{N+1}x_* - x_*\| \right)\..$$
      Therefore, at least one iteration of Landweber--Richardson method
      at level $N+1$ is performed and
      $$x_{k+1,N+1}^\delta = x_{k,N+1}^\delta - 
                             T_{N+1}^*
         \left( T_{N+1}(x_{k,N+1}^\delta)-y^\delta \right)\..$$
      The iteration at scale $X_{N+1}$ is terminated if 
      $$\|T_{N+1}(x_{k+1,N+1}^\delta) - y^\delta\| \leq 
        2(1+\eta)\left( \delta + \|P_{N+1}x_* - x_*\| \right)$$
      for the first time. The termination index 
      is denoted by $l_*(N+1,\delta)$.
\end{description}
\end{algorithm}
In the following, we verify some basic monotonicity properties of the 
iterates of the multi--level algorithm \ref{alg3.1}.
\begin{proposition}
\label{pr3.2}
Let Assumptions \ref{ass2.1CG} and \ref{ass2.2CG} hold.
At fixed level $N$ let the iteration be stopped 
according to the generalized discrepancy principle 
\be
\label{3.1}
\| \TN (x_{k,N}^\delta) - y^\delta\| > 
2(1+\eta) \left( \delta + \|P_Nx_*-x_*\| \right)\, \mbox{ for }
0 \leq k \leq l_*(N,\delta)\..
\ee
Then
\be
\label{3.2}
\|x_{j+1}^N - P_Nx_*\| \leq 
\|x_{j}^N - P_Nx_*\| \quad j=0,...,l_*(N,\delta)\;.
\ee
Moreover, if (\ref{3.1}) at level $N+1$ does not become 
active for $x_{0,N+1}^\delta$, then 
\be
\label{3.3}
\|x_{1,N+1}^\delta - x_*\|^2 \leq 
\|x_{0,N+1}^\delta - x_*\|^2 -  \frac{\eta}{1+\eta}
\|\TN x_{0,N+1}^\delta - y^\delta\|^2\..
\ee
\end{proposition}
{\bf{Proof:}}
From the definition of the multi--level Landweber--Richardson 
method it follows
\renewcommand{\arraystretch}{1.5}
\bea
\label{3.4}
\begin{array}{rcl}
\displaystyle
\|x_{k+1,N}^\delta-P_Nx_*\|^2 
&=& 
\displaystyle
\|x_{k,N}^\delta-P_Nx_*\|^2 + 
\|\TA(\TN(x_{k,N}^\delta)-y^\delta)\|^2 \\
& & 
\displaystyle
-2 \langle \TN(x_{k,N}^\delta-P_Nx_*),
          \TN(x_{k,N}^\delta)-y^\delta \rangle\\
&=& 
\displaystyle
\|x_{k,N}^\delta-P_Nx_*\|^2 - 2\| \TN(x_{k,N}^\delta)-y^\delta\|^2 \\
& & 
\displaystyle
+ \| \TA \left( \TN(x_{k,N}^\delta)-y^\delta \right)\|^2 \\
& & 
\displaystyle
+ 2 \langle \TN(P_Nx_*) - y^\delta, 
    \TN(x_{k,N}^\delta)-y^\delta \rangle \;.
\end{array}
\eea
From (\ref{3.1}) and (\ref{2.2CG}) it follows
\bea
\label{3.5}
\begin{array}{rcl}
\|\TN(P_Nx_*)-y^\delta\| 
&\leq&
\delta + \| T_{N_*}(P_Nx_*) - T_{N_*}(x_*) \| \\
&\leq&
\delta + \| P_Nx_* - x_*\|\..
\end{array}
\eea
For $0 \leq k \leq l_*(N,\delta)$ it follows from 
(\ref{2.1CG}), (\ref{3.4}) and (\ref{3.5}) that 
\bea
\label{3.6}
\begin{array}{rcl}
\displaystyle
\|x_{k+1,N}^\delta-P_Nx_*\|^2 
&\leq& 
\displaystyle
\|x_{k,N}^\delta-P_Nx_*\|^2 - \| \TN(x_{k,N}^\delta)-y^\delta\|^2 \\
& & 
\displaystyle
+ 2 \| \TN(P_Nx_*) - y^\delta \| 
    \| \TN(x_{k,N}^\delta)-y^\delta \|\\
&\leq& 
\displaystyle
\|x_{k,N}^\delta-P_Nx_*\|^2 - \frac{\eta}{1+\eta}
\| \TN(x_{k,N}^\delta)-y^\delta\|^2\;.
\end{array}
\eea
Now (\ref{3.2}) can be easily proven with an inductive argument.
(\ref{3.3}) can be proven with similar arguments as above.
~\hfill$\Box$
\bpn
From Proposition \ref{pr3.2} it follows that (\ref{3.1}) determines a 
well--defined stopping index for the Landweber--Richardson 
method at level $N$.
\begin{corollary}
\label{co3.3}
Let Assumptions \ref{ass2.1CG} and \ref{ass2.2CG} hold.
Let the Landweber--Richardson iteration at 
fixed level $N$ be stopped according to the generalized discrepancy 
principle (\ref{3.1}). 
If
\be
\label{3.7}
\delta + \| P_Nx_* - x_*\| > 0\,,
\ee
then 
$$ l_*(N,\delta) < \infty\;.$$
\end{corollary}
{\bf{Proof:}} From (\ref{3.1}) and (\ref{3.6}) it follows 
\beas
\begin{array}{rcl}
\displaystyle
4\eta(1+\eta) \sum_{k=0}^{l_*(N,\delta)} 
\left( \delta + \| P_Nx_* - x_*\| \right)^2
&\leq& 
\displaystyle
\|x_{l_*(N,\delta)+1,N}^\delta-P_Nx_*\|^2 \\
&    & 
\displaystyle
+ \frac{\eta}{1+\eta} \sum_{k=0}^{l_*(N,\delta)} 
\| \TN(x_{k,N}^\delta)-y^\delta\|^2 \\
&\leq& 
\displaystyle
\|x_{0,N}^\delta-P_Nx_*\|^2\\
&<& 
\displaystyle
\infty\,.
\end{array}
\eeas
This together with (\ref{3.7}) shows $l_*(N,\delta) < \infty$. 
\hfill$\Box$
\bpn
So far we have proven a monotonicity result for the iterates of the 
Landweber--Richardson method at fixed level $N$. 
In the following we show that the iterates of the multi--level 
algorithm also satisfy a certain monotonicity when the levels are 
switched. This together with Proposition \ref{pr3.2} provides 
monotonicity of the iterates of the multi--level algorithm~\ref{alg3.1}.
\bpn
\subsection{\label{sec3.1} Convergence and stability of the multi--level 
              algorithm}
Using the monotonicity results of the iterates of the multi--level 
schemes based on the Landweber--Richardson method we verify below that 
the multi--level scheme is convergent and stable.
The following lemma will be central for our further considerations. 
The proof of this lemma is analogous to the proof of Lemma 
\ref{le2.5CG}.
\begin{lemma}
\label{le3.4}
Let $\delta = 0$. 
If Assumptions \ref{ass2.1CG} and \ref{ass2.2CG} hold, the
multi--level algorithm \ref{alg3.1}
\begin{description}
\item{ } terminates after a finite number of iterations
\item{ or }
\item{ } 
      there exists a level $\hat{N}$ such that at each level
      $N < \hat{N}$ a finite number of Landweber--Richardson iterations
      are performed and at level $\hat{N}$ the 
      Landweber--Richardson iteration does not terminate.
\end{description} 
\end{lemma}
Similarily as in Section \ref{sec2} (using instead of Hanke's 
results in the proof of Theorem \ref{th2.6CG} some well--known 
results on convergence properties of the Landweber--Richardson 
method (see e.g. Groetsch \cite{Gro84}) it is possible to prove 
that the multi--level Landweber--Richardson algorithm is convergent 
in the case of exact data.
\begin{theorem}
\label{th3.5}
Let Assumptions \ref{ass2.1CG} and \ref{ass2.2CG} hold. Then the multi--level Landweber--Richardson 
algorithm converges to a solution $\hat{x}$. 
\end{theorem}
\bpn
In the case of perturbed data, it follows from Corollary \ref{co3.3} 
that the discrepancy principle (\ref{3.1}) determines a well--defined 
(finite) stopping index at level $N < \hat{N}$, which is denoted by 
$l_*(N,\delta)$. In the following we finally terminate the multi--level 
algorithm if for the first time
\be
\label{3.8}
\|\TN(x_{k,N}^\delta) - y^\delta\| \leq 
        2(1+\eta)\tau\delta\..
\ee
From Corollary \ref{co3.3} and Assumption \ref{ass2.2CG} it follows 
that this criterion determines a well--defined stopping index at level 
$N(\delta)$, which is denoted by $L_*(N(\delta),\delta)$.
\bpn
Our next result shows that the multi--level iteration is a stable method.
The proof of this result is analogous to the proof of Theorem 
\ref{th2.7CG}
\begin{theorem}
\label{th3.6}
Let the assumptions of Theorem \ref{th3.5} hold and assume that 
$T_{N_*}^*$ is continuously invertible. 
If $y^\delta$ fulfills (\ref{2.3CG}), and if the iteration at level $N$ 
is stopped according to the discrepancy principle (\ref{3.1}), and the 
multi--level iteration is terminated by (\ref{3.8}), then 
$$x_{N(\delta),L_*(N(\delta),\delta)}^\delta \to x_*, 
\qquad \delta \to 0\..$$
\end{theorem}  
\bpn
In practical applications, if the quantity $\|P^Nx_*-x_*\|$ 
is substantially underestimated, it may happen that 
the algorithm does not switch to the next level
although the iterations stagnate and the approximations do
not improve significantly any more. This can be avoided
for the multi--level iteration based on Landweber--Richardson 
method by using alternatively the following stopping criterion:
\be
\label{3.9}
\|x_{k+1}^N - x_k^N \| > 2(1+\eta) \tau \delta
  \mbox{ for } k=0,...,\hat{l}_*(N,\delta)\..
\ee
Then, from the definition of the Landweber--Richardson method 
it follows under the Assumption \ref{ass2.1CG}
$$\|(\TN x_N-y^\delta)\| \geq \|\TA(\TN x_N-y^\delta)\| 
  = \|x_{k+1}^N - x_k^N \| > 2(1+\eta) \tau \delta\..$$
Therefore we have $l_*(N,\delta) \geq \hat{l}_*(N,\delta)$. From 
this it is obvious to see that Proposition \ref{pr3.2},
Corollary \ref{co3.3}, Lemma \ref{le3.4} and Theorem \ref{th3.5} 
are valid, if we replace in these results $l_*$ by $\hat{l}_*$. 
It is also quite easy to see that (\ref{3.9}) determines a 
well--defined stopping criterion. The proof of Theorem \ref{th3.5}
for the multi--level Landweber--Richardson method with the 
stopping criterion \ref{3.9} requires the assumption that $T_{N_*}^*$ 
is continuously invertible to guarantee that the solution of 
$T_{N_*}^*x=y$ is a solution of $T_{N_*}(T_{N_*}^*-y)=0$ and vice versa.

\section{The multi--level algorithm for the moment problem} 
\label{sec4}
In the sequel we show the multi-level 
algorithms~\ref{alg2.3CG} and~\ref{alg3.1}
can be applied for the solution of general moment problems as 
described in Section \ref{sec1}. 
\bpn
It is well-known that a signal $x \in X_N$ can be completely
reconstructed from its moments $\langle x , g_j^N \rangle$ 
if the set $\{g_j^N\}_{j \in I}$ is a frame for $X_N$. For the 
readers convenience we quote the definition of a frame 
(see e.g. \cite{You80}):
\begin{definition}
\label{def2.7}
A family $\{g_j^N\}_{j \in I}$ is a called a frame for $X_N$ if 
there exist constants $A_N,B_N$, called the frame bounds, 
such that
\be
\label{18}
A_N\|x\|^2 \le \sum_{j \in I} |\langle x,g^N_j \rangle|^2
\le B_N \|x\|^2 \quad \forall x \in X_N \,.
\ee
\end{definition}
In this paper (if not stated otherwise) we will always assume 
that the family $\{g^N_j\}_{j \in I}$ is a frame for $X_N$.
For a given frame $\{g^N_j\}_{j\in I}$ we introduce the 
following linear operator
\be
\label{20}
\TN: X_N \rightarrow \ell^2(I)\,, \qquad
\TN x=\{\langle x, g^N_j \rangle\}_{j \in I}\..
\ee
The operator $\TN$ satisfies $\|\TN\| \le \sqrt{B_N}$. The 
adjoint operator is given by 
\be
\label{19}
\TA: \ell^2(I) \rightarrow X_N, \qquad 
\TA \{\mu_j\}_{j\in I}=\sum_{j \in I} \mu_j g^N_j \,.
\end{equation}
We note that if the family $\{g^N_j\}_{j \in J}$ is an orthonormal 
basis on $X_N$, then $\TN$ is a unitary operator and 
$\TA = (\TN)^{-1}$.
\bpn
As usual we define the {\it{frame operator}} by
\be
\label{21}
S_N: X_N \rightarrow X_N, \qquad S_Nx = \TA \TN  x =
 \sum_{j \in I}\langle x, g^N_j \rangle g^N_j \,.
\ee
It is well--known that $S_N$ is positive definite and invertible 
if $\{g_j^N\}_{j \in J}$ is a frame.
Stability results for frames in conjunction with moment
problems can be found in \cite{Chr96}.
\bpn

Obviously the operator $\TN$ (respectively the scaled operator
$\frac{1}{\sqrt{B}} T$) satisfies Assumption~\ref{2.1CG}(i) and (iii). 
A sufficient (but not necessary) condition for the 
frame analysis operators $\TN$ to satisfy Assumption~\ref{2.1CG}(ii)
is for instance when the frame elements $g^N_j$ correspond to reproducing 
kernels, see Section~\ref{sec4.1} below.
With the notation introduced above this shows the setting in which the
generally defined multi-level algorithms of Section~\ref{sec2} 
and Section~\ref{sec3} can be applied to recover $x$ from
noisy moment measurements
$y^{\delta}=\{\langle x,g_j^N \rangle +\delta_j \}_{j \in I}$.

For practical applications many delicate problems have
to be solved before implementation of the multi level algorithms.
We need conditions which 
allow an easy verification that the set $\{g_j^N\}$ is a frame
for $X_N$, as well as useful estimates of the frame bounds. And
finally we have to find appropriate estimates for $\|P_N x_* - x_*\|$
occurring in the stopping criteria~~(\ref{2.4CG}) and~(\ref{3.1}).

In Section 5 we summarize recent results from the literature
which allow to estimate frame bounds efficiently.

\subsection{Signal recovery in reproducing kernel Hilbert spaces} 
\label{sec4.1}
A special case of the moment problem occurs in reproducing kernel 
Hilbert spaces. 
Let $X$ be a Hilbert space of real or complex valued functions 
defined on a subset $I$ of $\r$ with inner product 
$\langle \cdot , \cdot \rangle$. 
We recall that a Hilbert space $X$ is called a 
{\it reproducing kernel Hilbert space} if for 
every element $t \in I$ the ``point evaluation'' 
$x \mapsto x(t)$ on $X$ is bounded. 
From the Riesz representation theorem it follows 
that there exists a unique $K_t \in X$, which satisfies
\be
\label{riesz}
\langle x,K_t \rangle = x(t) \quad \mbox{ for every } x \in X
\ee
The {\it{reproducing kernel}} $K$ on $I \times I$ is defined by
\be
\label{rk}
K(t,u) := \langle K_t, K_u \rangle = K_t(u)\..
\ee
In this section for application of the multi-level algorithm 
of Section \ref{sec2} we use a family of reproducing kernel 
Hilbert spaces $X_N$.
Note that every closed subspace of a RKHS is itself a RKHS.
Thus a reproducing kernel Hilbert space always contains a family
of nested reproducing kernel subspaces.

If $X_N \subseteq X$ is a RKHS, then the orthogonal projector 
of $X$ onto $X_N$ can be characterized by 
(see e.g. \cite[Chapter~11]{Ist87})

\be
\label{rkproj}
(P_{X_N} x)(t) = \langle x,K^N_t \rangle \..
\ee
The following lemma is an immediate consequence of the above 
characterization of the orthogonal projector $P_N$
\begin{lemma}
\label{le:proj}
Let $x_M \in X_M$, where $M \ge N$. Moreover, let $\TN$ be defined by 
$$ \TN(\cdot):= \langle \cdot, g_j^N \rangle\,,$$
where $g_j^N = K^N(\cdot,j)$. 
Then
\be
\label{proj}
\TN (P_N x_M) = \TN x_M\..
\ee
\end{lemma}
Using Lemma \ref{le:proj} we are able to verify the general 
assumptions in the convergence and stability results of the 
multi--level algorithm. Moreover, we discuss details of the 
numerical implementation, especially implementation of the 
stopping criteria. 

For reproducing kernel Hilbert space the stopping criterion 
(\ref{3.1}) is equivalent to 
$$
 |\sum_{j \in I}x_{k,N}^\delta (t_j) - y^\delta| \leq
      2(1+\eta) \left(\delta + \|P_Nx_* - x_*\|\right)\,.
$$

There are several important examples of reproducing kernel Hilbert spaces
which deserve to be mentioned explicitly.

(i)
Let $x \in {\cal B}_N$  where ${\cal B}_N$ is the 
{\it space of band-limited functions}, i.e.
\be
\label{bandlim}
\int_{-\infty}^{\infty} |x(t)|^2 dt < \infty \,\, \mbox{ and } \,\,
  \supp\, \hat{x} \subseteq [-N , N] \,,
\ee
where $\hat{x}$ denotes the Fourier transform of $x$.
The space ${\cal B}_N$ is a RKHS with reproducing kernel
\begin{equation}
K^{N}(t,u) = \frac{\sin N (t-u)}{N (t-u)} \,.
\end{equation}
This space is of paramount importance for digital signal processing,
see e.g.~\cite{Hig96}. We will deal with ${\cal B}_N$ in more detail in 
Section~\ref{sec5} in conjunction with signal recovery and bandwidth 
estimation from nonuniformly spaced noisy samples.

(ii)
The RKHS ${\cal P}_N$ consists of all {\it trigonometric 
polynomials}  $p$ of degree $N$ given by
\begin{equation}
p(t) = \sum_{n=-N}^{N} a_n e^{2 \pi i nt} \,.
\label{trigpol}
\end{equation}
The reproducing kernel is $K(t,u) = D_N(t-u)$ where the 
{\it Dirichlet kernel} is
\begin{equation}
D_N(t) = \frac{\sin \,(2N+1)\pi t}{\sin \, \pi t} \,.
\label{dirichlet}
\end{equation}
It follows immediately from the fundamental theorem of algebra that
the set $\{g_j\}_{j \in I}$, where $g_j(t) = D_N(t-t_j)$ constitutes a 
frame for ${\cal P}_N$ whenever $| I | > 2N$, where the sampling points
$t_j$ may be randomly spaced. Gr\"ochenig has derived estimates of the 
frame bounds for the weighted frame $\{w_j g_j\}_{j \in I}$,
if the maximal gap between two successive sampling points satisfies 
a certain criterion, see~\cite{Gro93}.

(iii)
{\it Sobolev spaces} $H^r$ for $r >\frac{1}{2}$ are well known to
possess reproducing kernels. Note however that the inner product
on $H^r$ is
not the restriction of that of $L^2(\r)$. Hence the formulae in the 
previous sections must be reinterpreted in terms of the $H^r$
inner product. See~\cite{NW91} for a detailed discussion of 
Sobolev spaces and reproducing kernels. 

{\it Wavelet subspaces} that constitute a so-called 
{\it multiresolution analysis} (MRA)~\cite{Mal89,Dau92} 
are further important examples for reproducing kernel Hilbert 
spaces~\cite{XiaNas94}, where moment problems arise
naturally when one wants to reconstruct a function from
its wavelet coefficients. Recently the concept of a MRA has 
been generalized, leading to a so-called {\it frame multiresolution 
analysis}~\cite{BL97}. 
More examples for reproducing kernel Hilbert spaces can be found in~\cite{NW91,Ist87}.

\section{Reconstructing a signal of unknown bandwidth from irregularly
spaced noisy samples}
\label{sec5}

In many practical applications, for instance in geophysics, spectroscopy, 
medical imaging, the signal cannot be sampled at regularly spaced points.
Thus one is confronted with the problem of reconstructing an
irregularly sampled signal. 
Often the signal can be assumed to be band-limited (or at
least essentially band-limited), but the actual bandwidth is not known.

The problem of reconstructing a band-limited signal $x$ from 
regularly and irregularly samples $\{x(t_j)\}_{j \in I}$ has 
attracted many mathematicians and engineers. Many theoretical results and 
efficient algorithms have been derived in the last years, 
see e.g.~\cite{Zay93,FG93,Sam95,Hig96} and the references cited 
therein. 
In \cite{Chr96} results on frame perturbations are applied to
the irregular sampling problem to derive error estimates for the 
approximation.
Most of these results are based on the assumption that the 
bandwidth of the signal to be recovered is known a priori.
However in many applications this assumption is not justified.
This is the motivation to consider the problem of reconstructing 
a band-limited signal whose actual band-width is not known, given
only a set of nonuniformly located and noisy samples.

In the sequel we show that the multi-level algorithm developed
in Section~\ref{sec2} in combination with the results of
Section~\ref{sec4} provides an efficient
solution to this problem.
To do this, we verify the general assumptions in the convergence and 
stability results of the multi--level algorithm.

In the approach suggested below, we use the
{\it{adaptive weights method}} developed by Feichtinger and 
Gr\"ochenig \cite{FG93}. The adaptive weights method is based on 
the fact that a signal $x \in \BN$ can be reconstructed from its 
nonuniformly spaced samples, if the maximal gap between two 
successive samples does not exceed the Nyquist interval $\pi/N$. The 
precise formulation of this results is summarized in the 
next theorem~\cite{FG93}.
\begin{theorem}
\label{th1}
Let $\{t_j\}_{j \in \z}$ be a strictly increasing sequence of real numbers 
(called the sampling set), which satisfies
\be
\label{maxgap}
\pi/N > \gamma := \sup_{j \in \z} (t_{j+1}-t_j)\,.
\ee
Then $\{\sqrt{w_j} K^{N}_{t_j}\}_{j \in \z}$ is a frame for
$B_{N}$ with frame bounds 
\begin{equation}
A = (1 - \frac{\gamma N}{\pi})^2 \mbox{ and } 
B = (1 + \frac{\gamma N}{\pi})^2\,.
\label{framab}
\end{equation}
where the weights $w_j$ are given by
\begin{equation}
w_j = (t_{j+1}-t_{j-1})/2 \,.
\end{equation}
\end{theorem}
\bpn
It is a consequence of Theorem \ref{th1} that a signal $x \in \BN$ 
can be recovered from its sampled data $\{x(t_j)\}$ by the 
Landweber--Richardson method, 
where the operators $\TN$ and $\TA$ are characterized by 
\be
\label{irr-sampl}
\TN x = \{\langle x,K^N_{t_j} \rangle \sqrt{w_j}\}_{j\in \z}\,\,,
\quad
\TA \{c_j\}_{j \in I} = \sum_{j \in I}  \sqrt{w_j} c_j K^N_{t_j} \,.
\ee
In the remaining of this section, we assume that for given sampling 
set $\{t_j\}_{j \in \z}$ the family $\{\sqrt{w_j} K_{t_j}^N\}_{j \in \z}$ 
satisfies (\ref{maxgap}) and thus forms a frame for $\BN$.

For an efficient implementation of the multi--level algorithm practical 
relevant estimates for $\|P_Nx_*-x_*\|$ are required. In this example
$P_N$ denotes the orthogonal projector of $L^2$ onto $\BN$.
Assume that the noisy samples $x_*^{\delta}(t_j) := x_*(t_j)+\delta_j$ are 
given with $\sum_{j \in \z} |\delta_j| = \delta$ and $x_* \in \BM$ for
fixed $M$. It holds
$$
\|x_* - P_N x_*\|^2 = \|x_*\|^2 - \|P_N x_*\|^2 = 
\int_{|\omega|>N}|\hat{x}_*(\omega)|^2 d\omega\.. 
$$
for $N \le M$, where we have used Plancherel's equality.
Here $\hat{x}$ denotes the Fourier transform of $x$.
In the sequel we show how to estimate $\|x_* - P_N x_*\|^2$ recursively.

It follows from the definition of frames that
\begin{equation}
\frac{1}{B_M} \sum_{j \in \z} |x_*(t_j)|^2 w_j \le \|x_*\|^2 \le
\frac{1}{A_M} \sum_{j \in \z} |x_*(t_j)|^2 w_j 
\label{framest1}
\end{equation}
with equality if the sampling points are regularly spaced.
In this case it is easy to see that 
$\{\sqrt{w_j} K_{t_j}^M\}_{j \in \z}$ is a tight frame for $\BN$
for $N \le M$.
Note that \ref{framest1} is only a worst-case estimate.
Extensive numerical experiments have shown that for practical 
applications
\begin{equation}
\label{framest2}
\sum_{j \in \z} |x^{\delta}_*(t_j)|^2 w_j - \|x_*\|^2 \approx \delta
\end{equation}
as long as $T$ is not very ill-conditioned.

If no a priori information about $\|P_N x_* - P_{N-1} x_*\|^2$ is available,
a reasonable estimate for $\|x_* - P_N x_*\|^2$ is given by
$\|x_* - P_{N-1} x_*\|^2$. 
Let us start with the first level, i.e. $N=1$. 
In this case we have 
$\|x_* - P_1 x_*\|^2  \le \|x_* - P_0 x_*\|^2 = \|x_*\|^2$.  
By virtue of~(\ref{framest2}) we can use
$\sum_{j \in \z} |x^{\delta}_*(t_j)|^2 w_j$ as estimate for 
$\|x_* - P_1 x_*\|^2$. Now we run algorithm~\ref{alg2.3CG}
or~\ref{alg3.1} at the first level until the corresponding stopping 
criterion~(\ref{2.4CG}) or~(\ref{3.1}) is satisfied.
We obtain an approximation denoted by $x_{k_*,1}^{\delta}$ and use 
$\|x_{k_*,1}^{\delta}\|^2$ as approximation to $\|P_2 x_*\|^2$
to estimate $\|x_* - P_2 x_*\|^2$ in the next level. We proceed 
inductively by using $\|x_{k_*,N}^{\delta}\|^2$ to estimate
\begin{equation}
\|x_* - P_{N+1} x_* \|^2 \approx 
\sum_{j \in \z} |x^{\delta}_*(t_j)|^2 w_j - \|x^{\delta}_{k_*,N}\|^2
\end{equation}
used in the stopping criterion at level $N+1$.

Improved estimates for $\|x_* - P_N x_*\|^2$ can be obtained
by trying to estimate $\|P_{N-1} x_* - P_N x_*\|^2 =
\int_{|\omega|=N-1}^{N} |\hat{x}_*(\omega)|^2 d\omega$ 
using a priori information about the decay of the frequencies
of $\hat{x}(\omega)$. This can be done either 
by a model-based approach or by a statistical analysis of a given class 
of signals. For instance potential fields, such as the gravity or 
the magnetic field, have exponentially decaying frequencies~\cite{RS97}. 
It is well-known that the Fourier coefficients of many bio-medical 
signals and large classes of images decay like $1/\omega^{\alpha}$ where
$\alpha$ is a constant depending on the class of signals under consideration.

\section{Experimental results} 
\label{sec6}

We demonstrate the performance of the proposed algorithms for the
reconstruction of a band-limited signal from its nonuniformly
spaced noisy samples. 
For the implementation of the algorithm we use the
discrete model presented in~\cite{FGS-act}.

\bpn
{\bf Experiment 1:}\enspace
In the first example we consider a signal from spectroscopy.
The original signal $x_*$ of bandwidth $M=30$ consists of 1024 
regularly spaced and nearly noise-free samples. We have added 
zero-mean white noise
with a signal-to-noise ratio of $12 \%$ and have sampled this noisy
signal at 107 nonuniformly spaced sampling points. 
A simple measure for the irregularity of the distribution
of the sampling points is the standard deviation $\sigma$ of the
differences of the sampling points. Let $\Delta_j = t_{j+1} - t_j$ 
and $s = \sum_{j=1}^{r} \Delta_j/r$ then $\sigma$ is given by
$$
\sigma = \sqrt{\frac{\sum_{j=1}^{r}|\Delta_j - s|^2}{r-1}} \,.
$$
Obviously $\sigma =0$ for regularly sampled signals.
In this experiment the maximal gap $\gamma$ in the sampling set is
33 and $\sigma =6.3$. The setup of this experiment 
simulates a typical situation in spectroscopy. 

Since we have to reconstruct $x_*$ from noisy samples, we cannot expect 
to recover $x_*$ completely. 

\begin{figure}
\begin{center}
\epsfig{file=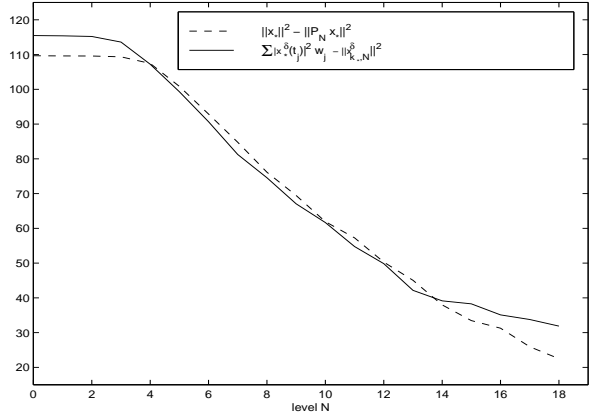,width=8cm,height=6cm}
\caption{Comparison of $\|x_* - P_N x_*\|^2$ and its estimation 
$\sum_{j} |x^{\delta}_*(t_j)|^2 w_j - 
\|x^{\delta}_{k_*,N}\|^2$.}
\label{fig:stopest}
\end{center}
\end{figure}

\begin{figure}
\begin{center}
\epsfig{file=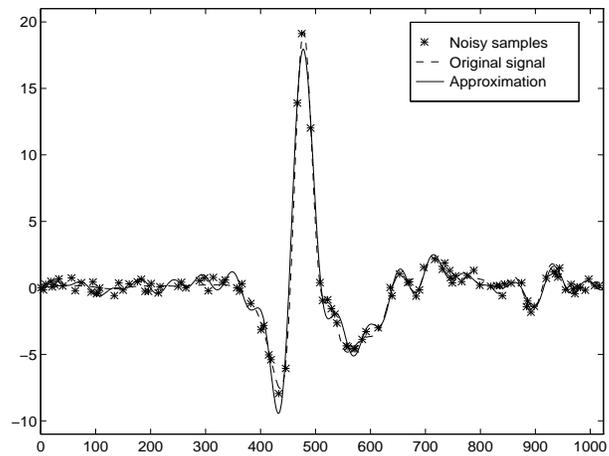,width=8cm,height=6cm}
\caption{Approximation using Algorithm~\ref{alg3.1}, 
         terminated at level 19.}
\label{fig:spec1}
\end{center}
\end{figure}

In our experiments we have not used any a priori information about the decay 
of the frequencies of the signal, but have estimated $\|x_* -P_N x_*\|$
by the simple recursive scheme as described in the previous section. 
A comparison of $\|x_* - P_N x_*\|^2$ and 
$\sum_{j} |x^{\delta}_*(t_j)|^2 w_j - \|x^{\delta}_{k_*,N}\|^2$
is demonstrated in Figure~\ref{fig:stopest}. 
Algorithm~\ref{alg3.1} terminates at level $N_*=19$, the
reconstruction is shown in Figure~\ref{fig:spec1} (solid line). 
The dashed line represents the original signal, and the $*$ are 
the noisy samples. The normalized error 
$\|x_* - x^\delta_{k_*,N}\|/\|x_*\|$ is $0.17$.

\begin{figure}
\begin{center}
\epsfig{file=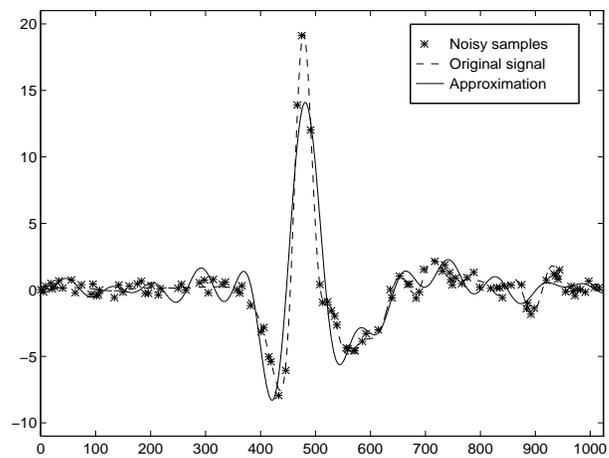,width=8cm,height=6cm}
\caption{Approximation using Algorithm~\ref{alg2.3CG}, terminating
         at level 13.}
\label{fig:spec2}
\end{center}
\end{figure}

\begin{figure}
\begin{center}
\epsfig{file=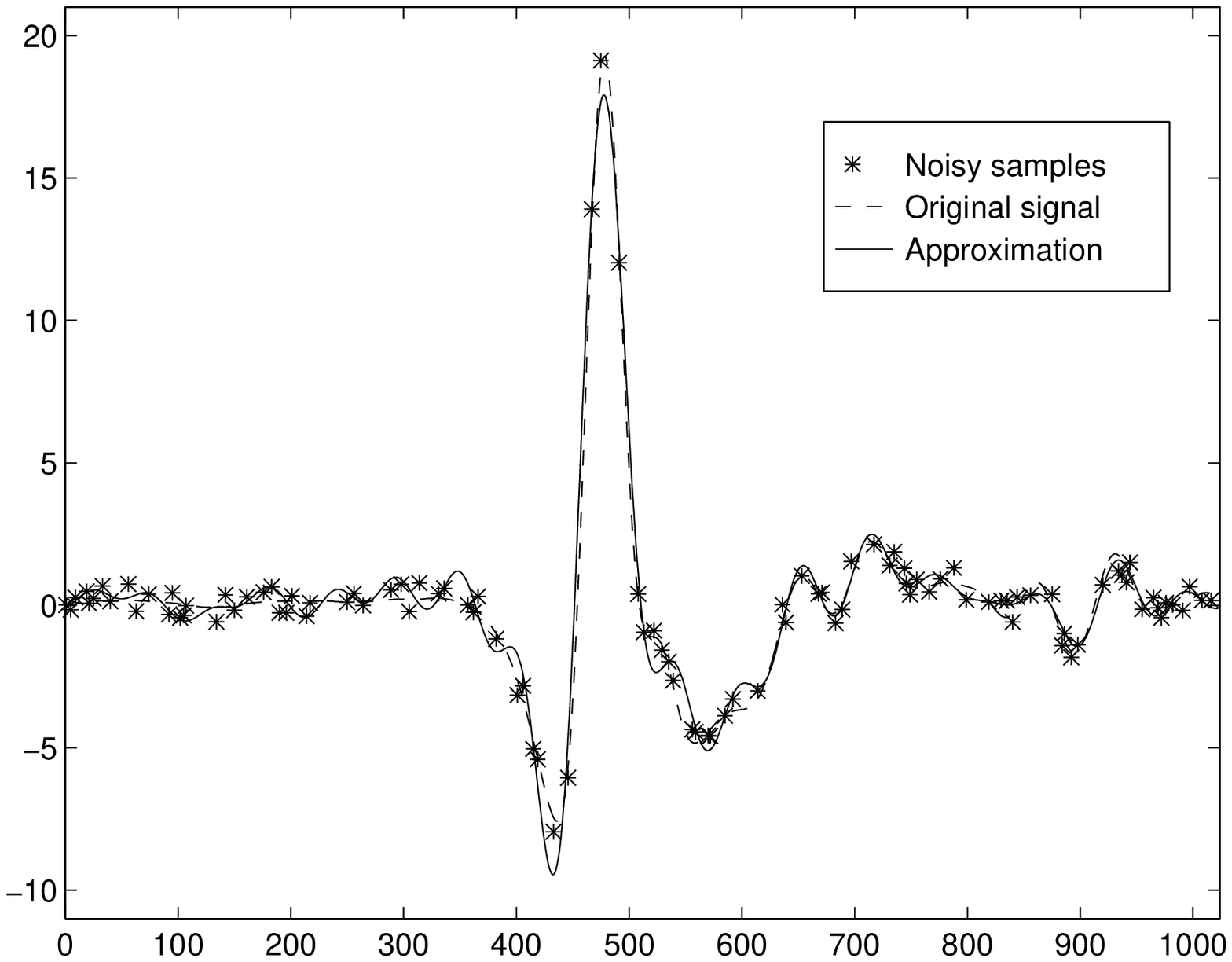,width=8cm,height=6cm}
\caption{Approximation using Algorithm~\ref{alg2.3CG}, using stopping
criterion~\ref{3.1}, terminating at level 19.}
\label{fig:spec3}
\end{center}
\end{figure}

Recall that we have had to use a different stopping criterion than 
for Algorithm~\ref{alg3.1} to prove convergence of the 
CG--multi-level algorithm. The CG-stopping criterion~(\ref{2.4CG}) 
forces the algorithm to terminate earlier,
thus giving a weaker approximation, see Figure~\ref{fig:spec2}.
Although we cannot justify it theoretically, we always got very good
reconstructions in our numerical experiments by using the same stopping criterion
as for the Landweber-Richardson multi-level algorithm, 
compare Figure~\ref{fig:spec3}.

\begin{figure}
\begin{center}
\epsfig{file=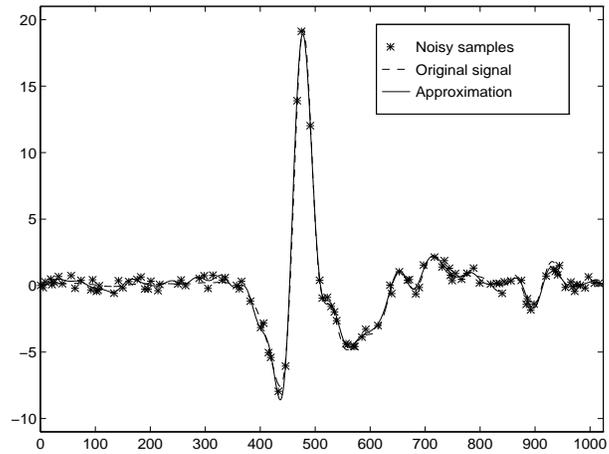,width=8cm,height=6cm}
\caption{Approximation using Algorithm~\ref{alg3.1}, omitting the 
factor 2 in the stopping criterion~\ref{3.1}, terminating at level 23.}
\label{fig:spec4}
\end{center}
\end{figure}

Note that our numerical experiments indicate that for most practical
applications one can drop the factor 2 in the discrepancy 
principle~(\ref{3.8}). In this case Algorithm~\ref{alg3.1} 
terminates at level 23 and the approximation error is $0.09$, 
see also Figure~\ref{fig:spec4}.

\begin{figure}
\begin{center}
\epsfig{file=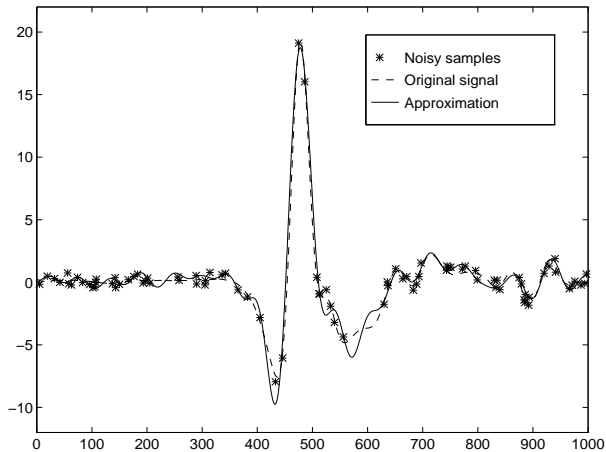,width=8cm,height=6cm}
\caption{Approximation using Algorithm~\ref{alg3.1} providing
an approximation error $\|x_* - x^{\delta}_{k_*,N_*}\|/\|x\| 
= 0.19$.}
\label{fig:spec5}
\end{center}
\end{figure}

\begin{figure}
\begin{center}
\epsfig{file=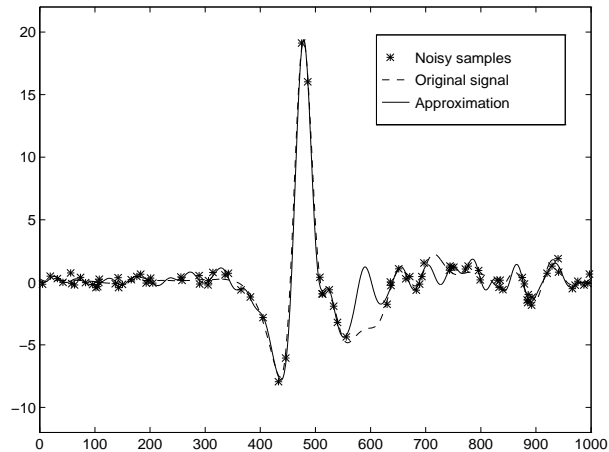,width=8cm,height=6cm}
\caption{Approximation using standard CGNE directly for level 
$M=30$, providing an approximation error 
$\|x_* - x^{\delta}_{N}\|/\|x\| = 0.28$.}
\label{fig:spec6}
\end{center}
\end{figure}

\bpn
{\bf Experiment 2:}\enspace
We use the same signal and the same noise level as in 
Experiment~1, but only 89 sampling points that are in addition 
more irregularly distributed with a standard deviation
$\sigma = 12.2$. Thus the operators
$\TN$ have a larger condition number than those arising
in Experiment~1. 
This experiment demonstrates that a proper choice
of the bandwidth can serve as a regularization procedure.
Algorithm~\ref{alg3.1} terminates at level $N_*=18$, the
final error between the original signal and the approximation
is 0.19. The reconstruction is shown in Figure~\ref{fig:spec5}.
It is remarkable that this approximation is better than the
approximation we obtain when we apply the standard CGNE
method~\ref{CGNE} directly using the a priori information
about the ``correct level'' $M=30$. In this case the
approximation error is 0.28, thus substantially larger.
The approximation is shown in Figure~\ref{fig:spec6}.

\bpn
{\bf Experiment 3:}\enspace
We consider the reconstruction of a signal representing a dynamical 
response of a mechanical system. The signal has been sampled at 300 
irregularly spaced points with $\sigma =40.7$, the original signal 
itself is not known. We run Algorithm~\ref{alg3.1} without using 
any other a priori information 
about the signal than the samples and the measurement error
which is known from experience to be approximately $10\%$. 

\begin{figure}
\begin{center}
\epsfig{file=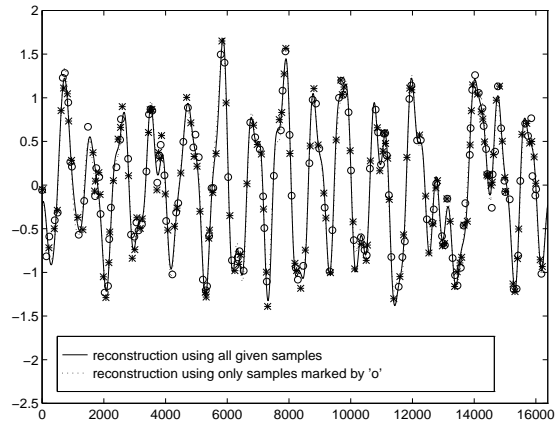,width=8cm,height=6cm}
\caption{Two reconstructions of dynamical impulse response 
using Algorithm~\ref{alg3.1}.}
\label{fig:dynamic}
\end{center}
\end{figure}

Algorithm \ref{alg3.1} terminates at level 42, the approximation 
is displayed in Figure~\ref{fig:dynamic} as solid line. In the 
chosen discrete model \cite{FGS-act}, $N=42$ means that the 
reconstruction
is a trigonometric polynomial of degree 42. Since the number of
sampling values is 300, we have about 3.5 times as many samples as we
would need theoretically to determine a reasonable approximation. 
Thus we should obtain about the same bandwidth and approximation by
using only every second sampling value, since we still have enough
samples and the corresponding operators $\TN$ are still well-conditioned.
It turns out that the algorithm terminates again at level 42,
giving an approximation close the first one.
The normalized error between both approximations is $11\%$.
Having in mind that we have used only $50\%$ of the given
samples and that the measurement error is about $10\%$, this
result demonstrates the robustness of the proposed algorithms.
The approximation for using only half of the samples is displayed
in Figure~\ref{fig:dynamic} as dotted line.

\section*{Acknowledgement}

Otmar Scherzer is partly supported by the Christian Doppler Society,
Austria, Thomas Strohmer is partly supported by the Austrian Science
Foundation FWF, Schr\"odinger scholarship J01388-MAT and project S7001-MAT.

We would like to thank S.~Spielman for providing the spectroscopy
data and Jacques Sainte-Marie for providing the dynamical impulse
response data. We also want to thank M.\ Hanke
(University of Kaiserslautern) for providing useful
information on CG methods.

\bibliographystyle{plain}

\end{document}